\documentclass{birkart}

\usepackage{amsmath}
\usepackage{amssymb}
\usepackage{amsfonts}

 \newtheorem{thm}{Theorem}[section]
 \newtheorem{cor}[thm]{Corollary}
 \newtheorem{lem}[thm]{Lemma}
 \newtheorem{prop}[thm]{Proposition}
 \theoremstyle{definition}

 \theoremstyle{remark}
 \newtheorem{rem}[thm]{Remark}

\numberwithin{equation}{section}
\numberwithin{figure}{section}

\newcommand{\cau}{{\mathfrak C}}
\newcommand{\beur}{{\mathfrak B}}
\newcommand{\proj}{{\mathfrak P}}
\newcommand{\integ}{{\mathfrak I}}
\newcommand{\Bergman}{A}
\newcommand{\ball}{\mathcal B}
\newcommand{\dbar}{\bar\partial}
\newcommand{\e}{\mathrm e}

\newcommand{\C}{{\mathbb C}}
\newcommand{\D}{{\mathbb D}}

\newcommand{\R}{{\mathbb R}}

\newcommand{\calL}{{X}}

\newcommand{\re}{\operatorname{Re}}
\newcommand{\im}{\operatorname{Im}}
\newcommand{\schlicht}{{\mathcal S}}

\newcommand{\diff}{{\mathrm d}}
\newcommand{\imag}{{\mathrm i}}

\newcommand{\const}{\Lambda}
\newcommand{\param}{\theta}
\newcommand{\maxop}{\mathfrak M}

\begin{document}
%
\title[Boundary properties of Green functions in the plane]
{Boundary properties of Green functions in the plane}
\author[Baranov]
{Anton Baranov}

\address{Baranov: Department of Mathematics and Mechanics\\
St. Petersburg State University\\
Universitetski\u\i{} prospekt, 28, \\
Stary\u\i{} Petergof, 198504 St. Petersburg\\
RUSSIA}
\email{a.baranov@ev13934.spb.edu}

\thanks{Research supported by the G\"oran Gustafsson Foundation.}

\author[Hedenmalm]
{H\aa{}kan Hedenmalm}

\address{Hedenmalm: Department of Mathematics\\
The Royal Institute of Technology\\
S -- 100 44 Stockholm\\
SWEDEN}

\email{haakanh@math.kth.se}




\begin{abstract}
We study the boundary properties of the Green function of bounded simply
connected domains in the plane. Essentially, this amounts to studying
the conformal mapping taking the unit disk onto the domain in question.
Our technique is inspired by a 1995 paper of Jones and Makarov.
The main tools are an integral identity as well as a uniform Sobolev 
imbedding theorem. The latter is in a sense dual to the exponential 
integrability of Marcinkiewicz-Zygmund integrals. 
We also develop a {\em Grunsky identity}, which contains the information of
the classical Grunsky inequality. This Grunsky identity is the case $p=2$
of a more general Grunsky identity for $L^p$ spaces.
\end{abstract}

\maketitle

\addtolength{\textheight}{2.2cm}







\section{Introduction}

\noindent{\bf Integral means spectrum of a conformal mapping.}
Let $\D$ denote the open unit disk in the complex plane $\C$. 
The class of univalent functions (conformal mappings) $\varphi:\D\to\C$,
subject to the normalizations $\varphi(0)=0$ and $\varphi'(0)=1$,
is denoted by $\schlicht$. Its subclass consisting of bounded
functions is denoted by $\schlicht_b$. 

For a given $\varphi\in\schlicht_b$ and a complex parameter $\tau$,
we consider the positive real numbers $\beta$ for which
\begin{equation}
\frac{1}{2\pi}\int_{-\pi}^{\pi}\big|\big[
\varphi'(r\e^{i\theta})
\big]^\tau\big| \,\diff\theta=O\bigg( \frac{1}{(1-r)^\beta}\bigg)
\qquad \text{as}\quad r\to 1^-,
\label{eq-intmean}
\end{equation}
holds.
Now, for a fixed $\tau$, we define the real number $\beta_\varphi(\tau)$ 
as the infimum of positive $\beta$ for which (\ref{eq-intmean}) holds.
In a sense, the quantity $\beta_\varphi(\tau)$ measures the 
expansion, compression, and rotation associated with the given conformal
mapping.
Finally, we define the {\it universal integral means spectrum} for the 
class $\schlicht_b$ by 
$$
\text{B}_b(\tau)=\sup_{\varphi \in \schlicht_b} \beta_\varphi(\tau).
$$
For some background information regarding the universal integral means 
spectrum, we refer to Makarov \cite{mak} as well as to Pommerenke's book
\cite{Pom}; see also Hedenmalm and Shimorin \cite{HedShi}, \cite{HedShiadd} 
as well as the survey paper by Hedenmalm and Sola \cite{HedSol}.
\medskip

\noindent \bf Weighted Bergman spaces. \rm For $-1<\alpha<+\infty$, 
we introduce the weighted Bergman spaces
${\Bergman}_\alpha^2(\D)$ consisting of those functions $f$, holomorphic in 
$\D$, for which
$$
\|f\|^2_{{\Bergman}_\alpha^2(\D)} =\int_\D |f(z)|^2 (1-|z|^2)^\alpha \,
\diff A(z)<+\infty;
$$
here, $\diff A(z)$ denotes normalized area measure in the plane,
$$\diff A(z)=\frac{\diff x\diff y}{\pi},\qquad z=x+\imag y.$$
The space ${\Bergman}_0^2(\D)$ obtained for parameter value $\alpha=0$ will
be denoted by ${\Bergman}^2(\D)$.
We point out the following almost obvious observation: 
\begin{equation}
\beta_\varphi(\tau) = \inf \Big\{\alpha+1:\, (\varphi')^{\tau/2}\in  
{\Bergman}_\alpha^2(\D) \Big\},
\qquad \tau\in\mathbb{C}.
\label{eq-beta}
\end{equation}

In Section \ref{alter}, we show that for $\varphi\in\schlicht_b$, 
\begin{equation}
\left|\log\frac{\varphi(z)}{z}\right|=O\left(\sqrt{\log\frac1{1-|z|^2}}
\right),\qquad |z|\to1^-,
\label{eq-loggrowth}
\end{equation}
so that, in particular, 
$$\left|\log\frac{\varphi(z)}{z}\right|=
o\left(\log\frac1{1-|z|^2}\right),\qquad |z|\to1^-.$$
In view of this, we see that 
\begin{equation*}
(\varphi')^{\tau/2}\in  
{\Bergman}_\alpha^2(\D)\quad\Rightarrow\quad
\bigg(\frac{z\varphi'}{\varphi}\bigg)^{\tau/2}\in{\Bergman}_\beta^2(\D)
\end{equation*}
holds whenever $\alpha<\beta$, while the opposite implication
\begin{equation*}
\bigg(\frac{z\varphi'}{\varphi}\bigg)^{\tau/2}\in{\Bergman}_\beta^2(\D)
\quad\Rightarrow\quad (\varphi')^{\tau/2}\in{\Bergman}_\alpha^2(\D)
\end{equation*}
holds whenever $\beta<\alpha$. The conclusion is that in terms of the 
integral means spectrum, it does not matter much if we consider the 
derivative $\varphi'$ or the function $z\varphi'(z)/\varphi(z)$. 

\medskip

\noindent\bf The estimate of the integral means spectrum. \rm
In this paper, two main results are obtained. The first runs as follows.

\begin{thm}
We have 
\begin{equation}
\text{\rm B}_b(2-\tau)\le 1 - \re \tau +
\bigg[\frac{9\e^2}{2}+o(1)\bigg]\,
|\tau|^2\log\frac1{|\tau|}\quad\text{as}\quad |\tau|\to0.
\end{equation}
\label{cor-main}
\end{thm}

Along the real line, Peter Jones and Nikolai Makarov
\cite{JonMak} obtained the better estimate
\begin{equation}
\text{\rm B}_b(2-\tau)\le1-t+O(t^2)\quad\text{as}\quad \R\ni t\to0;
\label{JM-est}
\end{equation}
in addition, they showed that for certain von Koch snowflake domains 
the $O(t^2)$ error term is sharp. These snowflakes should have fractal 
dimension arbitrarily close to $2$ as $t>0$ approaches $0$. As for Theorem
\ref{cor-main}, we thus expect von Koch-type snowflakes with spiraling 
behavior to be at least close to optimal. Other interesting candidates
for optimality are Fatou sets in iteration theory, as well as random domains
obtained by Loewner evolution, such as SLE (Schramm-Loewner evolution)
and DLA (diffusion limited aggregation). Standard SLE based on Brownian 
motion on the unit circle is known not to be optimal; nevertheless, other 
natural random motions on the unit circle still could do the job.
\medskip

\noindent{\bf The Grunsky inequality and operator identities.}
The classical Grunsky inequalities generalize Gr\"onwall's area theorem.
They constitute the most powerful tool known in the theory of univalent 
functions. In Section \ref{Grunsk}, we obtain an operator identity, which 
we have coined {\em Grunsky identity}; as we shall see, the Grunsky identity 
trivially implies the strong Grunsky inequality. 
The Grunsky identity involves 
Beurling-type singular integral operators. For $\varphi$ in the class 
$\schlicht$, let $\beur_\varphi$ be the singular integral operator
$$\beur_\varphi[f](z)=\text{pv}\int_\D
\frac{\varphi'(z)\varphi'(w)}{(\varphi(w)-\varphi(z))^2}\,f(w)\,\diff A(w),
\qquad z\in\D,$$
where ``pv'' indicates the standard principal value interpretation of the 
integral. From the fact that the standard Beurling operator
$$\beur_\C[f](z)=\text{pv}\int_\C\frac{f(w)}{(w-z)^2}\,\diff A(w),\qquad
z\in\C,$$
acts isometrically on $L^2(\C)$, we quickly deduce that $\beur_\varphi$
is a norm contraction on $L^2(\D)$. 
In the special case $\varphi(z)=z$, we write $\beur_e$ in place of 
$\beur_\varphi$. The Grunsky identity reads as follows:
\begin{equation}
\beur_\varphi-\beur_e=\proj\beur_\varphi=\beur_\varphi\bar\proj=
\proj\beur_\varphi\bar\proj,
\label{GI}
\end{equation}
where $\proj$ is the (Bergman) orthogonal projection on the analytic functions
in $L^2(\D)$, while $\bar\proj$ is the corresponding projection onto the
antianalytic functions. As the product of two contractions is a contraction,
we deduce from (\ref{GI}) that 
\begin{equation}
\big\|(\beur_\varphi-\beur_e)[f]\big\|_{L^2(\D)}\le\|f\|_{L^2(\D)},\qquad
f\in L^2(\D).
\label{GIn}
\end{equation} 
It is not too difficult to see that (\ref{GIn}) is an equivalent formulation
of the Grunsky inequalities (in the Grunsky inequalities, $f$ is assumed
antianalytic). For $0\le\theta\le2$, we consider the bigger collection of
Beurling-type operators 
$$\beur^\theta_\varphi[f](z)=\text{pv}\int_\D
\frac{\varphi'(z)^\theta\varphi'(w)^{2-\theta}}
{(\varphi(w)-\varphi(z))^2}\,f(w)\,\diff A(w),
\qquad z\in\D;$$
the case $\theta=1$ is the now familiar $\beur_\varphi$. For these operators,
we obtain the identity
\begin{equation}
\beur_\varphi^\theta-\beur_e+(\theta-1){\mathfrak M}_{1-|z|^2}
{\mathfrak M}_{\varphi''/\varphi'}\,{\mathfrak D}'=\beur_\varphi^\theta\bar
\proj,
\label{GI-theta}
\end{equation} 
where ${\mathfrak M}_F$ stands for the operator of multiplication by the 
function $F$, and the operator ${\mathfrak D}'$ is given by
$${\mathfrak D}'[f](z)=\int_\D\frac{f(w)}{(w-z)(1-\bar z w)}\,\diff A(w),
\qquad z\in\D.$$
The identity (\ref{GI-theta}) is the second main result of this paper.
For $1<p<+\infty$, let $K(p)$ denote the norm of $\beur_\C$ as an operator 
on $L^p(\C)$, and put $\theta=2/p$. It then follows from (\ref{GI-theta}) 
that
\begin{equation}
\big\|\big(\beur_\varphi^\theta-\beur_e+(\theta-1){\mathfrak M}_{1-|z|^2}
{\mathfrak M}_{\varphi''/\varphi'}\,{\mathfrak D}'\big)[f]\big\|_{L^p(\D)}
\le K(p)\,\|f\|_{L^p(\D)},
\label{GIn-p}
\end{equation} 
for $f$ in the subspace of $L^p(\D)$ consisting of the antianalytic functions.
The estimate (\ref{GIn-p}) should be thought of as an $L^p$ Grunsky 
inequality. 
\medskip

\noindent{\bf Reformulation in terms of the Green function.} 
Let $\Omega$ be a bounded simply connected domain in the complex plane $\C$. 
The phenomena considered in the
sequel are trivial for smooth boundaries $\partial\Omega$, so the interesting
thing is that the results apply to {\em all} bounded simply connected domains
$\Omega$, which may exhibit pretty wild boundary behavior. The Green function 
$G_\Omega(z,w)$ is, for a given $w\in\Omega$, the function which
is harmonic in $\Omega\setminus\{w\}$, vanishes along the boundary 
$\partial\Omega$ and has a logarithmic singularity of the form
$$G_\Omega(z,w)=\log(|z-w|^2)+O(1),$$ 
near $w$. By the strong maximum principle, the Green function is negative
at all points of $\Omega\times\Omega$ (with the understanding that it assumes 
the value $-\infty$ along the diagonal). Let us suppose the origin is 
contained in $\Omega$, and focus our attention to the function 
$G_\Omega(z)=G_\Omega(z,0)$. 
We introduce the standard Wirtinger differential operators
$$\partial_z=\frac12\bigg(\frac{\partial}{\partial x}-
\imag\frac{\partial}{\partial y}\bigg),\qquad
\bar\partial_z=\frac12\bigg(\frac{\partial}{\partial x}+
\imag\frac{\partial}{\partial y}\bigg),$$
as well as their multiplicative counterparts:
$$\partial_z^\times=z\partial_z,\qquad \bar\partial_z^\times=\bar z\bar 
\partial_z.$$
The function $\partial_z^\times G_\Omega(z)$ is analytic and zero-free in 
$\Omega$, and therefore has a well-defined analytic logarithm. Since 
$\partial_z^\times G_\Omega(z)|_{z=0}=1$, we may pick the logarithm that 
vanishes at $0$. We denote this logarithm by $H(z)$; thus
$$\e^{H(z)}=\partial_z^\times G_\Omega(z),\qquad H(0)=0.$$
We define the complex powers of $\partial_z^\times G_\Omega(z)$ in terms of 
the logarithm $H$:
$$\big[\partial_z^\times G_\Omega(z)\big]^\tau=\e^{\tau H(z)},\qquad 
\tau\in\C.$$
It is of interest to compare the sizes of 
\begin{equation}
\big|\big[\partial_z^\times G_\Omega(z)\big]^\tau\big|
\quad\text{and}\quad |G_\Omega(z)|^{-\alpha},
\label{eq-comp}
\end{equation}
where $\alpha$ is a real parameter. A natural way to do this 
is to use $L^1$ comparison: we say that given two positive Borel measurable 
functions $f,g$ on $\Omega$, $f$ is dominated by $g$ in $L^1(\Omega)$, in 
symbols, $f\ll_{L^1(\Omega)}g$, provided that $f/g\in L^1(\Omega)$. 
Moreover, the functions $f$ and $g$ are comparable in $L^1(\Omega)$, in 
symbols, $f\simeq_{L^1(\Omega)}g$, if $g\ll_{L^1(\Omega)}f\ll_{L^1(\Omega)}g$.
The question is for which real $\alpha$ and complex $\tau$ we have 
$$\big|\big[\partial_z^\times G_\Omega(z)\big]^\tau\big|
\ll_{L^1(\Omega)} |G_\Omega(z)|^{-\alpha},$$
that is,
\begin{equation}
\int_\Omega\big|\big[\partial_z^\times G_\Omega(z)\big]^\tau\big|\,
|G_\Omega(z)|^\alpha\,\diff A(z)<+\infty.
\label{eq-greenest}
\end{equation}
Let $\mathcal G_\Omega$ denote the subset of $\C\times\R$ consisting of pairs
$(\tau,\alpha)$ for which (\ref{eq-greenest}) holds. It is easy to see that
the set $\mathcal G_\Omega$ is convex, with the property that if 
$(\tau,\alpha)\in\mathcal G_\Omega$, then $(\tau,\alpha')\in\mathcal 
G_\Omega$ for all $\alpha'>\alpha$.
This means that there is a unique convex function ${\mathrm A}_\Omega:
\C\to\R$ such that if $\alpha<{\mathrm A}_\Omega(\tau)$, then 
$(\tau,\alpha)\notin\mathcal G_\Omega$, 
while for $\alpha>{\mathrm A}_\Omega(\tau)$, we have 
$(\tau,\alpha)\in\mathcal G_\Omega$. 

Let $\varphi:\D\to\Omega$ denote the conformal mapping which takes the open
unit disk onto $\Omega$, subject to the normalizations $\varphi(0)=0$, 
$\varphi'(0)>0$. The Green function is a conformal invariant, which means that
$$G_\Omega(\varphi(z))=\log(|z|^2),\qquad z\in\D,$$
and the change-of-variables formula for integrals shows that
\begin{equation}
\int_\Omega\big|\big[\partial_z^\times G_\Omega(z)\big]^\tau\big|\,
|G_\Omega(z)|^\alpha
\,\diff A(z)=\int_\D\bigg|\bigg[\frac{z\varphi'(z)}{\varphi(z)}\bigg]^{-\tau}
\bigg|\,\Big\{\log\frac1{|z|^2}\Big\}^\alpha|\varphi'(z)|^2\diff A(z).
\label{eq-greenid}
\end{equation}
We note that 
$$\log\frac1{|z|^2}=1-|z|^2+O\big((1-|z|^2)^2\big),\qquad |z|\to1^-,$$
which means that the right hand side of (\ref{eq-greenid}) essentially 
computes a weighted Bergman norm. 
In view of (\ref{eq-loggrowth}) and the remarks thereafter, we realize that
Theorem \ref{cor-main} is equivalent to the following.

\begin{cor}
The convex function ${\mathrm A}_\Omega(\tau)$ enjoys the following 
estimate:
$${\mathrm A}_\Omega(\tau)\le-\re\tau+\bigg[\frac{9\e^2}{2}+o(1)\bigg]\,
|\tau|^2\log\frac1{|\tau|}\quad\text{as}\quad |\tau|\to0,$$
where the $o(1)$ term is independent of the choice of the bounded simply
connected domain $\Omega$. 
\label{thm-main}
\end{cor}

If ${\mathrm A}_\Omega(\tau)+{\mathrm A}_\Omega(-\tau)\le0$, our scheme
of comparing the quantities in (\ref{eq-comp}) in terms of $L^1$ integrals
is very successful. It is therefore natural to view the quadratic-logarithmic
remainder term in Theorem \ref{thm-main} as the amount by which the $L^1$
comparison might fail. 

As a matter of fact, from a classical estimate of conformal mappings 
(\cite{Dur}, p. 126), we quickly derive that 
$${\mathrm A}_\Omega(\tau)\le\min\big\{|\tau-2|-1,\,|\tau|\big\},
\qquad\tau\in\C,$$
as well as 
\begin{equation}
{\mathrm A}_\Omega(\tau+\tau')\le{\mathrm A}_\Omega(\tau)+|\tau'|,\qquad
\tau,\tau'\in\C.
\label{eq-Aest}
\end{equation}
It is possible to use (\ref{eq-Aest}) to improve the estimate of Theorem 
\ref{thm-main} for $\tau$ with $\re \tau\le0$:
$${\mathrm A}_\Omega(\tau)\le-\re\tau+\bigg[\frac{9\e^2}{2}+o(1)\bigg]\,
|\im\tau|^2\log\frac1{|\im\tau|}\quad\text{as}\quad |\tau|\to0\quad\text{with}
\quad \re \tau\le0.$$
\medskip

\noindent{\bf Discussion of methods as regards Theorem \ref{cor-main}.} 
In the Jones-Makarov paper \cite{JonMak}, two methods are introduced. The
first (which seems to originate in the paper \cite{carl-jon}) is rather 
elementary, while the second uses very involved combinatorics, but yields the
strong estimate (\ref{JM-est}).
We generalize the first method only. It is based on an estimate of harmonic
measure ascribed to Arne Beurling, and a simple H\"older inequality estimate.
We replace the use of the estimate of harmonic measure by an identity 
involving the conformal mapping. The identity is the diagonal restriction of 
a more general identity which is an integrated version of the above-mentioned
Grunsky identity. 

It would be desirable to be able to obtain the smaller error term 
$O(|\tau|^2)$ in Corollary \ref{cor-main}. Quite possibly the combinatorial 
methods of Jones and Makarov might be adapted to achieve such a strengthening 
of the results. 
\medskip

\noindent{\bf Acknowledgements.} We gratefully acknowledge the generous
financial support of the G\"oran Gus\-tafs\-son foundation. We also thank Eero
Saksman for reading carefully a preliminary version of this manuscript.

\section{Transferred Cauchy and Beurling transforms}
\noindent{\bf The Cauchy transform and the Beurling transform.}
Let $\Omega$ be a bounded domain in the complex plane $\C$. 
For Lebesgue area integrable functions $f$ on $\Omega$,
we define the following two integral operators, the {\em Cauchy transform}
$\cau_\Omega$,
$$\cau_\Omega[f](z)=\int_\Omega \frac{f(w)}{w-z}\,\diff A(w),$$
and the {\em conjugate Cauchy transform} $\bar \cau_\Omega$
$$\bar \cau_\Omega[f](z)=\int_\Omega \frac{f(w)}{\bar w-\bar z}\,\diff A(w).$$
It is clear that in the sense of distribution theory,
$$\bar\partial_z \cau_\Omega[f](z)=-f(z),\qquad z\in\Omega,$$
and 
$$\partial_z \bar \cau_\Omega[f](z)=-f(z),\qquad z\in\Omega.$$

Associated to the Cauchy transform is the {\em Beurling transform}
$$\beur_\Omega[f](z)=\partial_z \cau_\Omega[f](z)=
\text{pv}\int_\Omega \frac{f(w)}{(w-z)^2}
\,\diff A(w),\qquad z\in\Omega,$$
while to the conjugate Cauchy transform we associate the {\em conjugate 
Beurling transform}
$$\bar \beur_\Omega[f](z)=\bar\partial_z \bar \cau_\Omega[f](z)=
\text{pv}\int_\Omega \frac{f(w)}{(\bar w-\bar z)^2}
\,\diff A(w),\qquad z\in\Omega.$$
It is well-known that for $\Omega=\C$, both $\beur_\C$ and $\bar \beur_\C$ 
are unitary transformations $L^2(\C)\to L^2(\C)$. Moreover, the adjoint 
$\beur_\C^*$ of $\beur_\C$ coincides with $\bar \beur_\C$, so that 
$\bar \beur_\C \beur_\C$
and $\beur_\C \bar \beur_\C$ equal the identity operator on $L^2(\C)$.
These assertions remain valid in case $\C\setminus\Omega$ has zero area.
For general $\Omega$, however, $\beur_\Omega$ and $\bar \beur_\Omega$ are 
contractions
$L^2(\Omega)\to L^2(\Omega)$, with $\beur_\Omega^*=\bar \beur_\Omega$.

Let $W^{1,2}(\Omega)$ be the space of all functions 
$f\in L^2(\Omega)$ such that
$\partial f\in L^2(\Omega)$ and 
$\bar\partial f\in L^2(\Omega)$ (the derivatives are understood 
in distributional sense).
We denote by  $W^{1,2}_{}(\Omega)/\C$ the 
Hilbert-Dirichlet space (modulo the constant functions) of functions $f$ on 
$\Omega$ with norm (the Dirichlet integral)
$$\|f\|_{W^{1,2}(\Omega)/\C}^2=\frac12\int_{\Omega}\big\{|\partial f(z)|^2+
|\bar\partial f(z)|^2\big\}\diff A(z).$$
Then $\cau_\Omega$ and $\bar \cau_\Omega$ are contractions 
$L^2(\Omega)\to W^{1,2}(\Omega)/\C$; in fact, this is an equivalent way of
saying that $\beur_\Omega$ and $\bar \beur_\Omega$ are contractions $L^2(\Omega)\to
L^2(\Omega)$. 

\begin{prop}
\label{prop-2.1}
For $f\in L^2(\Omega)$, extended to vanish on $\C\setminus\Omega$, we have 
$\cau_\C[f]\in W^{1,2}_{\text{\rm loc}}(\C)$.
\end{prop}

By the Sobolev inequality, functions in $W^{1,2}_{\text{loc}}(\C)$ belong to
all $L^{q}_{\text{loc}}(\C)$, $1<q<+\infty$. This may be improved 
substantially (see Chapter 3 of the book \cite{AH}).

\begin{prop}
For each $g\in W^{1,2}(\C)$ of norm $\le1$ in $W^{1,2}(\C)/\C$, we have 
$\exp[\beta_0|g|^2]\in L^1_{\text{\rm loc}}(\C)$, for some positive absolute 
constant $\beta_0$.
\label{prop-sobest}
\end{prop}

\noindent{\bf Conformal mapping and transferred Cauchy transforms.}
If $\Omega$ is simply connected (not the whole plane), there exists a 
conformal mapping $\varphi:\D\to\Omega$ that is onto. We connect two 
functions $f$ and $g$, on $\Omega$ and $\D$, respectively, via
$$g(z)=\bar\varphi'(z)\,f\circ\varphi(z),$$
and define the integral operator
$$\cau_\varphi[g](z)=(\cau_\Omega[f])\circ\varphi(z)=\int_\D 
\frac{\varphi'(w)}{\varphi(w)-\varphi(z)}\,g(w)\,\diff A(w),\qquad z\in\D;$$
$\cau_\varphi$ is then a contraction $L^2(\D)\to W^{1,2}(\D)/\C$, as 
follows from the previous observation that $\cau_\Omega$ is a contraction 
$L^2(\Omega)\to W^{1,2}(\Omega)$. Analogously, we define
$$\bar \cau_\varphi[g](z)=\int_\D 
\frac{\bar \varphi'(w)}{\bar \varphi(w)-\bar \varphi(z)}\,g(w)\,\diff A(w),
\qquad z\in\D,$$
and realize that $\bar \cau_\varphi$ is a contraction 
$L^2(\D)\to W^{1,2}(\D)/\C$. 
\medskip

\noindent\bf The transferred Beurling transforms. \rm 
It is well-known that $\beur_\C$ acts boundedly on $L^p(\C)$, for all $p$ with
$1<p<+\infty$. Let $K(p)$ be smallest positive constant 
\begin{equation}
\|\beur_\C f\|_{L^p(\C)}\le K(p)\,\| f\|_{L^p(\C)}, \qquad f\in L^p(\C).
\label{eq-Lp}
\end{equation}
The value of the optimal constant $K(p)$ is not known; however, it is known
that $K(2)=1$, and that in general 
$p^*-1\le K(p)\le2(p^*-1)$, where $p^*=\max\{p,p'\}$, and $p'=p/(p-1)$ is 
the dual exponent (see \cite{DV}). A conjecture of Tadeusz Iwaniec claims
that $K(p)=p^*-1$.

For $0\le\theta\le2$, we introduce the {\em $\theta$-skewed (transferred)
Beurling transform}, as defined by
\begin{equation}
\beur_\varphi^{\theta}[f]=
\text{pv}\int_\D\frac{\varphi'(z)^{\theta}\varphi'(w)^{2-\theta}}
{(\varphi(z)-\varphi(w))^2}\,f(w)\,\diff A(w).
\label{eq-skewedb}
\end{equation}
It follows from (\ref{eq-Lp}) that
$$\big\|\beur_\varphi^{2/p}[f]\big\|_{L^p(\D)}\le K(p)\,\|f\|_{L^p(\D)},\qquad
f\in L^p(\D),$$
for all $p$ with $1<p<+\infty$. In the symmetric case $\theta=1$, we shall 
write $\beur_\varphi$ in place of $\beur_\varphi^1$. We note that 
$\beur_\varphi$ is
a contraction on $L^2(\D)$.

\section{An operator identity related to the Grunsky inequality}
\label{Grunsk}

\noindent\bf The basic identity. \rm 
In case $\varphi(z)=z$, we write $\cau_e$ instead of $\cau_\varphi$, and 
$\bar \cau_e$ in place of $\bar \cau_\varphi$. Likewise, under the same 
circumstances,
we write $\beur_e$ and $\bar \beur_e$ instead of $\beur_\varphi$ and 
$\bar \beur_\varphi$.
Next, let $\proj:L^2(\D)\to A^2(\D)$ be the orthogonal projection,
$$\proj [f](z)=\int_\D \frac{f(w)}{(1-\bar w z)^2}\,\diff A(w),$$
and let $\integ_0$ be the operation on analytic functions in $\D$ defined by
$$\integ_0[f](z)=\int_0^z f(w)\,\diff w,\qquad z\in\D.$$
Analogously, we let $\bar\proj$ be the orthogonal projection to the 
antianalytic functions in $L^2(\D)$, and we let $\bar\integ_0$ be the 
corresponding integration operator acting on the antianalytic functions. 
In terms of formulas, we have
$$\bar\proj [f](z)=\int_\D \frac{f(w)}{(1-w \bar z)^2}\,\diff A(w),$$
and
$$\bar\integ_0[f](z)=\int_0^z f(w)\,\diff \bar w,\qquad z\in\D.$$

The following identity is basic to our investigation.

\begin{prop}
For $\varphi\in\schlicht$, we have the identity
\begin{equation}
\log\frac{z(\varphi(z)-\varphi(\zeta))}{(z-\zeta)\varphi(z)}+
\log(1-\bar z\zeta)=\int_\D
\frac{\varphi'(w)}{\varphi(w)-\varphi(z)}\,\frac{\zeta}{1-\bar w \zeta}\,
\diff A(w).
\label{eq-fund}
\end{equation}
\label{prop-1}
\end{prop} 

\begin{proof}
We note that for analytic functions $f$ area-integrable in $\D$, we have 
$$\int_\D f(w)\,\frac{\zeta}{1-\bar w\zeta}\,
\diff A(w)=\int_0^\zeta f(w)\,\diff w.$$
It follows that
$$\int_\D\left[
\frac{\varphi'(w)}{\varphi(w)-\varphi(z)}-\frac{1}{w-z}\right]\,
\frac{\zeta}{1-\bar w \zeta}\,\diff A(w)
=\log\frac{z(\varphi(z)-\varphi(\zeta))}{(z-\zeta)\varphi(z)}.$$
Next, we compute that
$$\int_\D\frac{\zeta}{(w-z)(1-\bar w \zeta)}\,\diff A(w)=
\log\big(1-\bar z\zeta\big).$$
The assertion is now immediate.
\end{proof}

\medskip

\noindent\bf The Grunsky identity. \rm
If we apply the differentiation operator $\partial_\zeta$ to the identity of 
Proposition \ref{prop-1}, we get
$$\frac{\varphi'(\zeta)}{\varphi(\zeta)-\varphi(z)}-
\frac{1}{\zeta-z}
-\frac{\bar z}{1-\bar z\zeta}=\int_\D
\frac{\varphi'(w)}{\varphi(w)-\varphi(z)}
\,\frac{1}{(1-\bar w\zeta)^2}\,\diff A(w),$$
which in terms of operators may be written in the form 
\begin{equation}
\cau_\varphi-\cau_e-\bar \integ_0\bar \proj=\cau_\varphi \bar\proj.
\label{eq-fundop}
\end{equation}
As we apply the operator $\partial$ (differentiation with respect to $z$) 
to both sides, we obtain the derived identity
\begin{equation}
\beur_\varphi-\beur_e=\beur_\varphi \bar\proj.
\label{eq-fundop-2}
\end{equation}
We call (\ref{eq-fundop-2}) a {\em Grunsky identity}, and claim that it
{\em trivially entails the strong Grunsky inequality}. 
The main observation needed is that
both $\beur_\varphi$ and $\bar\proj$ are contractions on $L^2(\D)$, making 
their product $\beur_\varphi\bar\proj$ a contraction as well. If $f$ is in 
$L^2(\D)$, then
\begin{equation}
\big\|(\beur_\varphi-\beur_e)[f]\big\|^2_{L^2(\D)}=
\big\|\beur_\varphi\bar\proj [f]\big\|^2_{L^2(\D)}\le\|f\big\|^2_{L^2(\D)}.
\label{eq-fundop-3}
\end{equation}
It is an easy exercise to check that (\ref{eq-fundop-3}) -- applied to an
antianalytic $f$ -- is equivalent to the strong Grunsky inequality as 
formulated in \cite{Dur}. To help the reader, we indicate how this is done. 
First, there is the step of translating between the classes $\schlicht$
and $\Sigma$ of conformal mappings, which is standard. Next, let
$$\frac{\varphi'(z)\varphi'(w)}{(\varphi(w)-\varphi(z))^2}-\frac{1}{(w-z)^2}
=\sum_{j,k=0}^{+\infty}\gamma_{j,k}\,z^jw^k$$
be the usual Taylor series expansion of two variables of the left hand side; 
in terms of the Grunsky coefficients $\beta_{j,k}$ as presented in 
\cite{Dur}, we have 
$$\gamma_{j,k}=-(k+1)\,\beta_{j+1,k+1}.$$
We apply (\ref{eq-fundop-3}) to the antianalytic polynomial
$$f(z)=\sum_{j=0}^{N}a_j\,\bar z^j.$$
This results in
\begin{equation}
\sum_{j=0}^{+\infty}\frac1{j+1}\bigg|\sum_{k=0}^{N}
\frac{\gamma_{j,k}}{k+1}\,a_k\bigg|^2\le\sum_{j=0}^{N}
\frac{|a_j|^2}{j+1}.
\label{GI-D}
\end{equation}
The estimate (\ref{GI-D}) now expresses the classical formulation of the 
Grunsky inequalities. 
\medskip

\noindent\bf A skewed Grunsky identity. \rm
A formulation of the Grunsky identity which is pretty much equivalent to 
(\ref{eq-fundop-2}) runs as follows:
\begin{equation}
\beur_\varphi-\beur_e=\proj\beur_\varphi.
\label{eq-fundop-2'}
\end{equation}
Let us try to find an analog of (\ref{eq-fundop-2'}) for the $\theta$-skewed 
Beurling transform $\beur_\varphi^\theta$, as defined by (\ref{eq-skewedb}),
with $0<\theta<2$. 
First, we note that
$$\frac{\varphi'(z)^{\theta}\varphi'(w)^{2-\theta}}
{(\varphi(z)-\varphi(w))^2}=\frac{1}{(z-w)^2}+(\theta-1)
\frac{\varphi''(w)}{(z-w)\varphi'(w)}+O(1)$$
holds near the diagonal $z=w$, which we interpret to say that
$$\beur_\varphi^\theta-\beur_e+(\theta-1)
\cau{\mathfrak M}_{\varphi''/\varphi'}$$
maps into the analytic functions. Here, ${\mathfrak M}_{\varphi''/\varphi'}$
is the operator of multiplication by the function $\varphi''/\varphi'$.
In other words, we can show that 
\begin{equation}
\beur_\varphi^\theta-\beur_e+(\theta-1)
\cau{\mathfrak M}_{\varphi''/\varphi'}=\proj\big(
\beur_\varphi^\theta-\beur_e+(\theta-1)
\cau{\mathfrak M}_{\varphi''/\varphi'}\big)=\proj\beur_\varphi^\theta+
(\theta-1)\proj\cau{\mathfrak M}_{\varphi''/\varphi'}.
\label{eq-prel}
\end{equation}
We readily calculate that
$\proj\cau=\integ_0^*$, where $\integ_0^*$ is the integral operator
$$\integ_0^*[f](z)=\int_\D\frac{\bar w}{1-\bar w z}\,f(w)\,\diff A(w),$$
which is in a reasonable sense the adjoint of the integration operator
$\integ_0$. As a consequence, we find that
$$(\cau-\integ_0^*)[f](z)=\int_\D\frac{1-|w|^2}{(w-z)(1-\bar w z)}\,f(w)\,
\diff A(w),\qquad z\in\D.$$
We factor $\cau-\integ_0^*={\mathfrak D}{\mathfrak M}_{1-|z|^2}$, where
${\mathfrak M}_{1-|z|^2}$ is the operator of multiplication by $1-|z|^2$, and
$\mathfrak D$ is defined by
$${\mathfrak D}[f](z)=\int_\D\frac{f(w)}{(w-z)(1-\bar w z)}\,
\diff A(w),\qquad z\in\D.$$
{\em We finally obtain from (\ref{eq-prel}) the skewed Grunsky identity}
\begin{equation}
\beur_\varphi^\theta-\beur_e+(\theta-1)
{\mathfrak D}{\mathfrak M}_{1-|z|^2}{\mathfrak M}_{\varphi''/\varphi'}=
\proj\beur_\varphi^\theta.
\label{skGrunsky}
\end{equation}
The operators $\proj$, $\beur_e$, and $\beur^{2/p}_\varphi$ are all bounded 
on $L^p(\D)$, for $1<p<+\infty$. It can be shown that $\mathfrak D$
is also a bounded operator on $L^p(\D)$, for $1<p<+\infty$. 
As a matter of fact, it is possible to read this off from (\ref{skGrunsky}) 
(we refrain from supplying the necessary details). 
The special case $\theta=1$ of (\ref{skGrunsky}) is indeed
(\ref{eq-fundop-2'}). 

A variant of (\ref{skGrunsky}) reads as follows:
$$\beur_\varphi^\theta-\beur_e+(\theta-1){\mathfrak M}_{1-|z|^2}
{\mathfrak M}_{\varphi''/\varphi'}\,{\mathfrak D}'=\beur_\varphi^\theta\bar
\proj,$$
where ${\mathfrak D}'$ is given by
$${\mathfrak D}'[f](z)=\int_\D\frac{f(w)}{(w-z)(1-\bar z w)}\,\diff A(w),
\qquad z\in\D.$$
As a sample (skewed) Grunsky-type estimate which may be obtained in this 
fashion, we mention
\begin{multline}
\int_\D\bigg|\frac{\varphi'(z)^{2/p}\varphi'(w)^{2-2/p}}
{(\varphi(z)-\varphi(w))^2}-\frac1{(z-w)^2}+\bigg(\frac2p-1\bigg)
\,\frac{\varphi''(z)(1-|z|^2)}{(w-z)(1-\bar z w)\varphi'(z)}\bigg|^p\diff A(z)
\\
\le K(p)^p\,{}_2F_1\big(p,p;2;|w|^2\big),\qquad w\in\D,
\label{eq-exGr}
\end{multline}
where $K(p)$ is as in (\ref{eq-Lp}), and ${}_2F_1$ is Gauss' hypergeometric
function. The case $p=2$ of (\ref{eq-exGr}) is an invariant version of 
Gr\"onwall's classical area theorem \cite{Gron}. It is an easy exercise to
derive from (\ref{eq-exGr}) the pointwise estimate
\begin{multline*}
\bigg|\frac{\varphi'(z)^{2/p}\varphi'(w)^{2-2/p}}
{(\varphi(z)-\varphi(w))^2}-\frac1{(z-w)^2}+\bigg(\frac2p-1\bigg)
\,\frac{1}{(z-w)^2}\log\frac{\varphi'(w)}{\varphi'(z)}\bigg|
\\
\le K(p)\,\big\{{}_2F_1\big(p,p;2;|w|^2\big)\big\}^{1/p}
\big\{{}_2F_1\big(p',p';2;|z|^2\big)\big\}^{1/p'},\qquad z,w\in\D,
\end{multline*} 
where $p'=p/(p-1)$ is the dual exponent.
\medskip

\noindent\bf Variants of the basic identity. \rm
If we restrict Proposition \ref{prop-1} to the diagonal $z=\zeta$, we get 
the following.

\begin{cor}
For $\varphi\in\schlicht$, we have the identity
\begin{equation*}
\log\frac{z\varphi'(z)}{\varphi(z)}+
\log\big(1-|z|^2\big)=\int_\D
\frac{\varphi'(w)}{\varphi(w)-\varphi(z)}\,\frac{z}{1-\bar w z}\,
\diff A(w).
\end{equation*}
\label{cor-1}
\end{cor} 

For the applications we have in mind, it will be convenient to work with
the following variant of Proposition \ref{prop-1}. 

\begin{prop}
We have the identity
\begin{multline*}
\log\frac{z(\varphi(z)-\varphi(\zeta))}{(z-\zeta)\varphi(z)}
-\zeta(1-|\zeta|^2)\bigg[\frac{\varphi'(\zeta)}{\varphi(\zeta)-\varphi(z)}
-\frac{1}{\zeta-z}\bigg]+\log\big(1-\bar z\zeta\big)+
\bar z\zeta\,\frac{1-|\zeta|^2}{1-\bar z\zeta}\\
=\zeta^2\int_\D
\frac{\varphi'(w)}{\varphi(w)-\varphi(z)}\,\frac{\bar\zeta-\bar w}
{(1-\bar w\zeta)^2}\,\diff A(w).
\end{multline*}
\label{prop-2}
\end{prop}

\begin{proof}
We have that
$$\zeta^2\,\frac{\bar \zeta-\bar w}{(1-\bar w \zeta)^2}=
\frac{\zeta}{1-\bar w \zeta}-
\frac{\zeta(1-|\zeta|^2)}{(1-\bar w \zeta)^2},$$
so that for analytic functions $f$ area-integrable in $\D$, we get
$$\zeta^2\int_\D f(w)\,\frac{\bar \zeta-\bar w}{(1-\bar w \zeta)^2}\,
\diff A(w)=\int_0^\zeta f(w)\,\diff w-\zeta(1-|\zeta|^2)f(\zeta).$$
In particular, this leads to the identity
\begin{multline*}
\zeta^2\int_\D\left[
\frac{\varphi'(w)}{\varphi(w)-\varphi(z)}-\frac{1}{w-z}\right]\,
\frac{\bar \zeta-\bar w}{(1-\bar w \zeta)^2}\,\diff A(w)\\=
\log\frac{z(\varphi(z)-\varphi(\zeta))}{(z-\zeta)\varphi(z)}
-\zeta(1-|\zeta|^2)\bigg[\frac{\varphi'(\zeta)}{\varphi(\zeta)-\varphi(z)}
-\frac{1}{\zeta-z}\bigg].
\end{multline*}
Next, we compute that
$$\zeta^2\int_\D\frac{\bar\zeta-\bar w}{(w-z)(1-\bar w \zeta)^2}\,\diff A(w)=
\log\big(1-\bar z\zeta\big)+\bar z\zeta\,\frac{1-|\zeta|^2}{1-\bar z\zeta}.$$
The assertion is now immediate.
\end{proof}

By plugging in the diagonal $z=\zeta$ in Proposition \ref{prop-2}, we
arrive at the following.
 
\begin{cor}
We have the identity
\begin{multline*}
\log\frac{z\varphi'(z)}{\varphi(z)}
-z(1-|z|^2)\,\frac{\varphi''(z)}{\varphi'(z)}+\log\big(1-|z|^2\big)+|z|^2
=z^2\int_\D
\frac{\varphi'(w)}{\varphi(w)-\varphi(z)}\,\frac{\bar z-\bar w}
{(1-\bar w z)^2}\,\diff A(w).
\end{multline*}
\label{cor-2}
\end{cor}

We show that the additional term which appears in the identity
of Proposition \ref{prop-2} is uniformly bounded with respect to $z$
and $\zeta$.

\begin{lem}
There is an absolute constant $C>0$ such that
\begin{equation}
(1-|\zeta|^2)\bigg|\frac{\varphi'(\zeta)}{\varphi(\zeta)-\varphi(z)}
-\frac{1}{\zeta-z}\bigg| \le C
\end{equation}
for any $z,\zeta\in\D$.
\label{lem-3.3}
\end{lem}

\begin{proof} By a classical property of the class $\schlicht$,
for any $\psi\in\schlicht$, we have
$$
\bigg|w\frac{\psi'(w)}{\psi(w)}\bigg|\le\frac{1+|w|}{1-|w|}, \qquad w\in\D,
$$
and, consequently,
\begin{equation}
\label{lem3.3a}
\bigg|\frac{\psi'(w)}{\psi(w)}-\frac{1}{w}\bigg|\le\frac{C_1}{1-|w|^2}, 
\qquad w\in\D,
\end{equation}
for some absolute constant $C_1$. We apply (\ref{lem3.3a}) to the function
$$
\psi(w) = \frac{\varphi\big(\frac{w+z}{1+\bar z w}\big)-\varphi(z)}
{(1-|z|^2)\varphi'(z)},
$$
and after the change of variables 
$$\zeta=\frac{w+z}{1+\bar z w},$$
we obtain
$$
\bigg|\frac{\varphi'(\zeta)}{\varphi(\zeta)-\varphi(z)}
-\frac{1}{\zeta-z}\frac{1-|z|^2}{1-\zeta \bar z}\bigg| 
\le\frac{C_1}{1-|\zeta|^2},\qquad z,\zeta\in \D.
$$
This is clearly equivalent to the asserted estimate.
\end{proof}

\section{The growth of the argument for bounded univalent functions}
\label{alter}

\noindent{\bf An estimate of the logarithm of $\varphi(z)/z$.}
It is clear that for $\varphi\in\schlicht_b$, the analytic function
\begin{equation*}
\Phi(z)=\log\frac{\varphi(z)}{z},\qquad z\in\D,
\end{equation*}
has bounded real part. However, it is easy to see from examples that the 
imaginary part, which corresponds to taking the argument of $\varphi(z)/z$, 
can be unbounded. It is therefore a natural question to ask for growth bounds
of the function $\Phi$. To this end, we consider the derivative
$$\Phi'(z)=\frac{\varphi'(z)}{\varphi(z)}-\frac1{z},\qquad z\in\D.$$
As $\varphi$ is bounded and univalent, $\varphi'$ has bounded $L^2(\D)$ 
integral, that is, $\varphi'\in A^2(\D)$. Since in the annulus 
$\frac12<|z|<1$, $\varphi(z)$ is bounded from below by the standard 
estimates, we get that $\varphi'(z)/\varphi(z)$ is $L^2$-integrable in the 
annulus. Since $\Phi(z)$ is analytic throughout $\D$, it follows that 
$\Phi'\in A^2(\D)$. By the reproducing property of the Bergman kernel, we have
$$\Phi'(z)=\int_\D\frac{\Phi'(w)}{(1-\bar w z)^2}\,\diff A(w),\qquad z\in\D,$$
and since $\Phi(0)=0$, we may integrate both sides to obtain
$$\Phi(z)=z\int_\D\frac{\Phi'(w)}{1-\bar w z}\,\diff A(w),\qquad z\in\D.$$
By the Cauchy-Schwarz inequality, we get
\begin{equation*}
|\Phi(z)|\le|z|\,\|\Phi'\|_{A^2(\D)}\left\{\int_\D\frac{\diff A(w)}
{|1-\bar w z|^2}\right\}^{1/2}=\|\Phi'\|_{A^2(\D)}
\sqrt{\log\frac{1}{1-|z|^2}},\qquad z\in\D.
\end{equation*}
In particular, we see that
\begin{equation*}
\bigg|\log\frac{\varphi(z)}{z}\bigg|=
O\bigg(\sqrt{\log\frac{1}{1-|z|^2}}\bigg),\qquad |z|\to1^-.
\end{equation*}

\section{Marcinkiewicz-Zygmund integrals}
\label{MZ}

\noindent \bf Zygmund's paper. \rm
Here, we basically recall the estimates of Antoni Zygmund \cite{Zyg}, where
he mentions that a part of the results was inspired by a remark by his
former student Richard O'Neil.
Fix a bounded domain $\Omega$ in $\C$, and let $\kappa$, $0<\kappa<+\infty$,
be a real parameter. For $w\in\C$, let $\delta(w)$ denote the Euclidean 
distance from $w$ to $\C\setminus\Omega$; for $w\in\C\setminus\Omega$, then,
$\delta(w)=0$. Pick a real parameter $\param$, confined to the interval 
$0<\param<1$.
The Marcinkiewicz-Zygmund integral is defined by the formula
$$I_\kappa(z)=\int_\Omega\min\bigg\{\frac{\delta(w)^\kappa}{|z-w|^{2+\kappa}},
\frac{\param^{-2-\kappa}}{\delta(w)^{2}}\bigg\}\,\diff A(w).$$
Let $g$ be a positive locally area-integrable function in $\C$. The {\em
Hardy-Littlewood maximal function} for $g$ is defined by
$$\maxop[g](z)=\sup_{0<r<+\infty}\frac1{r^2}\int_{\D(z,r)}g(w)\,\diff A(w),$$
where $\D(z,r)$ stands for the open disk of radius $r$ with center point $z$.
We need the associated function
$$h(r,w)=\int_{\D(w,r)}g(z)\,\diff A(z);$$
we clearly have
$$h(r,w)\le r^2\,\maxop[g](w).$$ 
By an integration by parts argument, we have that for $0<\rho<+\infty$,
\begin{equation*}
\frac1{\rho^2}\int_{\D(w,\rho)}g(z)\,\diff A(z)+
\rho^\kappa\int_{\C\setminus\D(w,\rho)}\frac{g(z)}{|z-w|^{2+\kappa}}\,
\diff A(z)
=(2+\kappa)\rho^\kappa\int_\rho^{+\infty}h(r,w)\,r^{-3-\kappa}\diff r.
\end{equation*}
By comparing with the maximal function, we get
\begin{equation}
\frac1{\rho^2}\int_{\D(w,\rho)}g(z)\,\diff A(z)+
\rho^\kappa\int_{\C\setminus\D(w,\rho)}\frac{g(z)}{|z-w|^{2+\kappa}}\,
\diff A(z)
\le\frac{2+\kappa}{\kappa}\,\maxop[g](w).
\label{eq-max}
\end{equation}
By Fubini's theorem,
\begin{multline*}
\int_\C I_\kappa(z)\,g(z)\,\diff A(z)=
\int_\Omega\bigg\{\frac{\param^{-2-\kappa}}{\delta(w)^2}
\int_{\D(w,\param\delta(w))}g(z)\,
\diff A(z)\\
+\delta(w)^\kappa\int_{\C\setminus\D(w,\param\delta(w))}
\frac{g(z)}{|z-w|^{2+\kappa}}\,\diff A(z)\bigg\}\diff A(w),
\end{multline*}
so that in view of (\ref{eq-max}), we obtain
\begin{equation}
\int_\C I_\kappa(z)\,g(z)\,\diff A(z)
\le \param^{-\kappa}\,\frac{2+\kappa}{\kappa}\int_\Omega \maxop[g](w)\,
\diff A(w).
\label{eq-max1}
\end{equation}
Let $L(p)$ be a constant for which
\begin{equation}
\|\maxop[g]\|_{L^{p'}(\C)}\le L(p)\,\|g\|_{L^{p'}(\C)},\qquad g\in L^{p'}(\C),
\label{maxest}
\end{equation}
where $p$ and $p'$ are dual exponents: $p'=p/(p-1)$. In Grafakos' book 
\cite{graf}, it is shown that
$$L(p)=3^{2/p'}\,p\le 9p$$
works. As we combine (\ref{eq-max1}) with the maximal function estimate 
(\ref{maxest}), we arrive at
\begin{equation*}
\int_\C I_\kappa(z)\,g(z)\,\diff A(z)
\le \param^{-\kappa}\frac{2+\kappa}{\kappa}\,|\Omega|_A^{1/p}\,
\|\maxop[g]\|_{L^{p'}(\C)}
\le \param^{-\kappa}
\frac{2+\kappa}{\kappa}\,L(p)\,|\Omega|_A^{1/p}\,\|g\|_{L^{p'}(\C)},
\end{equation*}
which amounts to an $L^p$ estimate of $I_\kappa$:
\begin{equation*}
\|I_\kappa\|_{L^p(\C)}\le \param^{-\kappa}\frac{2+\kappa}{\kappa}\,L(p)\,
|\Omega|_A^{1/p}.
\end{equation*}
Here, $|\Omega|_A$ denotes the normalized area of $\Omega$.
We pick a complex parameter $\lambda$, and estimate
\begin{equation}
\big\|\e^{\lambda I_\kappa}-1\big\|_{L^1(\C)}\le
\sum_{j=1}^{+\infty}\frac{1}{j!}\,|\lambda|^j\,\|I_\kappa\|_{L^j(\C)}^j\le
|\Omega|_A\sum_{j=1}^{+\infty}
\frac1{j!}\,|\lambda|^j\,L(j)^j\,\param^{-j\kappa}
\bigg[\frac{2+\kappa}{\kappa}\bigg]^j.
\label{exp-est}
\end{equation}
The bound from Grafakos' book means that the $L^1$ norm of 
$\e^{\lambda I_\kappa}-1$ is bounded in a well-controlled manner for 
$$|\lambda|<\frac{\kappa\,\param^\kappa}{9\e(2+\kappa)}.$$
Indeed, a crude version of Stirling's formula yields
\begin{equation} 
\big\|\e^{\lambda I_\kappa}-1\big\|_{L^1(\C)}\le
\frac{\kappa\,|\Omega|_A}{\kappa-9\e|\lambda|\param^{-\kappa}(2+\kappa)}
-|\Omega|_A
\quad\text{for}\quad|\lambda|<\frac{\kappa\,\param^\kappa}{9\e(2+\kappa)}.
\label{exp-est1}
\end{equation}

\section{Uniform Sobolev imbedding}

\noindent\bf A modified transferred Cauchy transform. \rm
We introduce a relative of the transferred Cauchy transform $\cau_\varphi$,
as defined by 
$$
\widetilde{\cau}_\varphi [f](z)
= \int_\D \frac{\varphi'(w)}{\varphi(w)-\varphi(z)}\,\frac{\bar z-\bar w}
{1-\bar w z} f(w)\,\diff A(w).
$$
This operator is related to the right hand side expression in Corollary 
\ref{cor-2}. It also solves the $\dbar^2$ problem
$$\dbar^2\,\widetilde{\cau}_\varphi [f](z)=-\frac{f(z)}{1-|z|^2},\qquad
z\in\D.$$
\medskip

\noindent\bf Application of H\"older's inequality. \rm
For positive $\kappa$, we consider the Lebesgue space
$$\calL_\kappa(\D) = L^p(\D,\mu),$$
where 
$$p=\frac{2+\kappa}{1+\kappa},\qquad 
\diff\mu(z)=(1-|z|^2)^{-\kappa/(1+\kappa)}\,\diff A(z);$$
the norm in the Banach space $\calL_\kappa(\D)$ is given by
$$
\|f\|_{\calL_\kappa(\D)} =
\Bigg\{ \int_\D |f(z)|^{(2+\kappa)/(1+\kappa)}
(1-|z|^2)^{-\kappa/(1+\kappa)}\, \diff A(z)
\Bigg\}^{(1+\kappa)/(2+\kappa)}.
$$

Suppose that $f\in \calL_\kappa(\D)$. Then, by H\"older's inequality,
\begin{equation}
\label{marc1}
\big| \widetilde{\cau}_\varphi [f](z) \big|\le \bigg\{\int_\D
\bigg|\frac{(w-z)\varphi'(w)}{(1-\bar w z)(\varphi(w)-\varphi(z))}
\bigg|^{2+\kappa}(1-|w|^2)^\kappa\, \diff A(w) \bigg\}^{1/(2+\kappa)}\\
\times \big\| f \big\|_{\calL_\kappa(\D)}.
\end{equation}
\medskip

\noindent\bf Marcinkiewicz-Zygmund integrals. \rm
The function
$$J_{\kappa}[\varphi](z)=\int_\D
\bigg|\frac{(w-z)\varphi'(w)}{(1-\bar w z)(\varphi(w)-\varphi(z))}
\bigg|^{2+\kappa}(1-|w|^2)^\kappa
\,\diff A(w)$$
is essentially the familiar Marcinkiewicz-Zygmund integral.
Indeed, put $\Omega=\varphi(\D)$, and consider the 
Marcin\-kie\-wicz-Zyg\-mund integral $I_\kappa$ associated to $\Omega$
(see the previous section).
Recall that $\delta(z)$ denotes the Euclidean distance from $z\in\Omega$ to
the complement $\C\setminus\Omega$. It is well known that 
$$\frac14\,(1-|z|^2)|\varphi'(z)|\le\delta(\varphi(z))\le
(1-|z|^2)|\varphi'(z)|,\qquad z\in\D.$$
In the integral defining $J_\kappa[\varphi](z)$, a suitably large hyperbolic 
disk about $z$ may be deleted from the area of integration, as it contributes 
only a bounded quantity to the integral. Changing the variables, then, it 
is now easy to check that
$$J_k[\varphi](z)\le4^{\kappa}\,I_\kappa(\varphi(z))
+O_{\param,\kappa}(1),\qquad z\in\D.$$
In view of Section \ref{MZ}, we know that for bounded $\Omega$,
$$
\int_\Omega  \e^{|\lambda|\, I_\kappa(z)}\, \diff A(z)\le
\frac{\kappa\,|\Omega|_A}{\kappa-9\e|\lambda|\param^{-\kappa}(2+\kappa)}
\quad\text{for}\quad|\lambda|<\frac{\kappa\,\param^\kappa}{9\e(2+\kappa)}.
$$
This leads to an integral estimate of $J_\kappa[\varphi]$ on $\D$:
\begin{equation}
\int_{\D}\e^{|\lambda|\,J_{\kappa}[\varphi](z)} 
\,|\varphi'(z)|^2\, \diff A(z) \le C(\param,\kappa)\|\varphi'\|^2_{L^2(\D)}
\frac{\kappa}{\kappa-9\e|\lambda|4^\kappa \param^{-\kappa}(2+\kappa)},
\label{eq-expint-2}
\end{equation}
for 
$$|\lambda|<\frac{\kappa 4^{-\kappa}\param^\kappa}{9\e(2+\kappa)},$$
where $C(\param,\kappa)$ denotes an appropriate positive constant that 
depends only on $\param$ and $\kappa$.
\medskip

\noindent\bf A uniform Sobolev imbedding theorem. \rm
Denote by ${\ball}_\kappa$ the unit ball in the space $\calL_\kappa(\D)$.
It follows from (\ref{marc1}) that
$$
\sup\limits_{f\in {\ball}_\kappa} \big|\widetilde{\cau}_\varphi 
[f](z)\big|^{2+\kappa}
\le J_k[\varphi](z), \qquad z\in\D.
$$
As the parameter $\param$ may be chosen arbitrarily close to $1$, we have 
-- in view of (\ref{eq-expint-2}) -- obtained the following estimate.

\begin{lem}
For any positive $\kappa$, we have, for complex $\lambda$ with
$$|\lambda|<\frac{\kappa 4^{-\kappa}}{9\e(2+\kappa)},$$
the estimate
\begin{equation*}
\int_{\D}  \exp\Big\{|\lambda|\sup\limits_{f\in {\ball}_\kappa}
\big|\,\widetilde{\cau}_\varphi[f](z)\,\big|^{2+\kappa} \Big\}
\,|\varphi'(z)|^2\, \diff A(z)<+\infty.
\end{equation*}
\label{lem-sobtype}
\end{lem}
\smallskip

\begin{rem}
The inequality of Lemma \ref{lem-sobtype} should be compared with the 
following corollary of the classical Sobolev inequality (see Proposition
\ref{prop-sobest}).
For $\zeta\in\D$, consider the function  
$g_\zeta(z) = \zeta/ (1- \bar z \zeta)$. 
Note that the right-hand side of our basic identity (Proposition 
\ref{prop-1}) is $\cau_\varphi[g_\zeta](z)$.
It is easy to see that
$$
\|g_\zeta\|_{L^2(\D)}^2 = \log\frac{1}{1-|\zeta|^2}.
$$
Now put $f_\zeta =g_\zeta/\|g_\zeta\|_{L^2(\D)}$.
Then, by Propositions \ref{prop-2.1} and \ref{prop-sobest} (carried out on
the target domain $\varphi(\D)$), we have
\begin{multline}
\int_\D \exp\Big\{\beta_0\, 
\big|\cau_\varphi[f_\zeta](z)\big|^2\Big\}\,
|\varphi'(z)|^2\, \diff A(z)  
=\int_\D \exp\Bigg\{\beta_0\frac{|
\cau_\varphi[g_\zeta](z)|^2}{\log\frac{1}{1-|\zeta|^2}}
\Bigg\}\,  |\varphi'(z)|^2\, \diff A(z) <+\infty.
\label{sob2}
\end{multline}
If we could establish a ``diagonal'' analog of inequality (\ref{sob2})
with $\zeta=z$ (so that $\zeta$ is not constant in the integral), this  would 
imply the strong Jones-Makarov estimate (with $O(t^2)$ in the error term). 
However, the only substitute we have is the uniform Sobolev imbedding of 
Lemma \ref{lem-sobtype} which indeed allows us to plug in the diagonal 
choice $\zeta=z$, but at the price of weakening the estimate. See the next 
section for details.
\end{rem}
\medskip

\section{The proof of the main theorem}

\noindent{\bf Application of Lemma \ref{lem-sobtype}.}
\rm
We now specialize to $z=\zeta$ in Proposition \ref{prop-2} (see also Corollary
\ref{cor-2}). In view of Lemma \ref{lem-3.3}, we obtain 
\begin{equation}
\label{special}
\log\frac{z\varphi'(z)}{\varphi(z)} + \log\big(1-|z|^2) +O(1)
=z^2\int_\D
\frac{\varphi'(w)}{\varphi(w)-\varphi(z)}\,\frac{\bar z-\bar w}
{(1-\bar w z)^2}\,\diff A(w), \qquad z\in\D.
\end{equation}
Note that the right-hand side of (\ref{special}) is 
$\widetilde\cau_\varphi [g_z](z)$, with $g_z(w)= z^2/(1-\bar w z)$. 
Thus, we have
\begin{equation}
\label{upper}
\widetilde\cau_\varphi [g_z](z) =
\log\Bigg[ z\frac{\varphi'(z)}{\varphi(z)}\,\big(1-|z|^2\big)\Bigg]
+O(1), \qquad z\in\D.
\end{equation}
We plan to apply Lemma \ref{lem-sobtype} to 
the function $f=f_z=g_z/\|g_z\|_{\calL_\kappa(\D)}$. 
As a first step, we estimate the norm
$$
\|g_z\|_{\calL_\kappa(\D)}^{2+\kappa} =
|z|^{2(1+\kappa)} \Bigg(\int_\D\frac{(1-|w|^2)^{-\kappa/(1+\kappa)}}
{|1-\bar w z|^{(2+\kappa)/(1+\kappa)}}\,
\diff A(w)\Bigg)^{1+\kappa}
$$
for $z\in\D$ with $|z|$ close to $1$. Note that for fixed $\vartheta$ with 
$-\frac12<\vartheta<+\infty$,
\begin{equation*}
\int_\D\frac{(1-|w|^2)^{2\vartheta}}{|1-\bar w z|^{2+2\vartheta}}\,\diff A(w)
=(1+o(1))\,\frac{\Gamma(1+2\vartheta)}{\Gamma(1+\vartheta)^2}\,
\log\frac1{1-|z|^2},
\quad\text{as}\quad|z|\to1^-.
\end{equation*}
We now permanently restrict $\kappa$ to the interval $0<\kappa<1$, and
conclude that whenever the positive constant $\const$ is chosen so that
$$\const>\bigg[\frac{\Gamma\big(\frac{1-\kappa}{1+\kappa}\big)}
{\Gamma\big(\frac1{1+\kappa}\big)^2}\bigg]^{1+\kappa},$$
we obtain that
\begin{equation}
\label{lower}
\|g_z\|_{\calL_\kappa(\D)}^{2+\kappa} \le \const\,
\Bigg[\log\frac1{1-|z|^2}\Bigg]^{1+\kappa}+O(1), \qquad z\in\D;
\end{equation}
here, we assume that $\kappa$ is restricted to the interval $0<\kappa<1$
(after all, we will be only interested in small $\kappa$). 
By combining (\ref{upper}) with (\ref{lower}), we obtain
\begin{equation*} 
\bigg| \log\bigg[ z\frac{\varphi'(z)}{\varphi(z)}(1-|z|^2)\bigg]
\bigg|^{2+\kappa}
\bigg/ \bigg[\log\frac{1}{1-|z|^2}\bigg]^{1+\kappa}
\le \const\,\big| \widetilde\cau_\varphi[f_z](z) 
\big|^{2+\kappa}
+O(1),\qquad z\in\D.
\end{equation*}
The function $f_z\in\ball_\kappa$ is as before given by 
$f_z=g_z/\|g_z\|_{\calL_\kappa(\D)}$.
We now apply Lemma \ref{lem-sobtype} to obtain, for $\lambda\in\C$ with
$$|\lambda|<\frac{\kappa 4^{-\kappa}}{9\e(2+\kappa)},$$
\begin{multline}
\label{est1}
\int_{\D}
\exp \Bigg\{\frac{|\lambda|}{\const}\,\Big| 
\log\Big[ z\frac{\varphi'(z)}{\varphi(z)}(1-|z|^2)\Big]\Big|^{2+\kappa}
\Big/ \Big[\log\frac{1}{1-|z|^2}\Big]^{1+\kappa} \Bigg\}\, |\varphi'(z)|^2
\diff A(z)
\\
= \int_{\D}
\exp \Bigg\{\frac{|\lambda|}{\const} \,
\Bigg| 1- \frac{\log\frac{z\,\varphi'(z)}{\varphi(z)}}
{\log\frac{1}{1-|z|^2}} \Bigg|^{2+\kappa}\,\log\frac{1}{1-|z|^2} 
\Bigg\}\, |\varphi'(z)|^2\diff A(z)
<+\infty.
\end{multline}

\noindent\bf Linear approximation argument. \rm
We use a very simple argument to complete the proof of Theorem 1.1.
We apply the convexity estimate
\begin{multline*}
|a|^{2+\kappa}=|\bar a|^{2+\kappa}\ge |b|^{2+\kappa}-(2+\kappa)|b|^\kappa
\re\big[b(\bar b-a)\big] \\
=|b|^{2+\kappa}+(2+\kappa)|b|^\kappa
\big[\re b-|b|^2\big]-(2+\kappa)|b|^\kappa
\re\big[ b(1-a)\big], \qquad a,b\in\mathbb{C},
\end{multline*}
to 
$$
a= 1- \frac{\log\frac{z\varphi'(z)}{\varphi(z)}}
{\log\frac{1}{1-|z|^2}},
$$
and obtain
\begin{multline}
\Bigg| 1- \frac{\log\frac{z\,\varphi'(z)}{\varphi(z)}}
{\log\frac{1}{1-|z|^2}} \Bigg|^{2+\kappa}\,\log\frac{1}{1-|z|^2} 
\\
\ge
\Big[|b|^{2+\kappa}+(2+\kappa)|b|^\kappa
\big[\re b-|b|^2\big]\Big]\log\frac1{1-|z|^2}
-(2+\kappa)|b|^\kappa
\re\bigg[b \log\frac{z\varphi'(z)}{\varphi(z)}\bigg]
\label{eq-linest}
\end{multline}
for any $b\in\mathbb{C}$. As we insert the estimate (\ref{eq-linest}) into
(\ref{est1}), we find that
\begin{multline*}
\int_{\D}
\exp \Bigg\{ \frac{|\lambda|}{\const} 
\Big[|b|^{2+\kappa}+(2+\kappa)|b|^\kappa
\big[\re b-|b|^2\big]\Big]\log\frac1{1-|z|^2} \\
-
\frac{|\lambda|}{\const}\,(2+\kappa)\,|b|^\kappa\,\re\bigg[b 
\log\frac{z\varphi'(z)}{\varphi(z)}\bigg]
\Bigg\}
\,|\varphi'(z)|^2 \diff A(z) <+\infty.
\end{multline*}
Next, we assume $b\neq0$, and put 
$\tau = \const^{-1} |\lambda|(2+\kappa)\,|b|^\kappa b$.
Note also that
$$
\exp \Bigg\{-\frac{|\lambda|}{\const}\,(2+\kappa)\,|b|^\kappa\,
\re\bigg[b\log\frac{z\varphi'(z)}{\varphi(z)}
\bigg]\Bigg\} = \bigg| \Big[\frac{z\varphi'(z)}{\varphi(z)}\Big]^{-\tau} 
\bigg|.
$$
Thus, we have 
\begin{equation}
\label{est2}
\int_{\D}
\bigg| \Big[ \frac{z\varphi'(z)}{\varphi(z)}\Big]^{-\tau} \bigg|\,
\big(1-|z|^2\big)^{-\re \tau +M\,|\tau|^{(2+\kappa)/(1+\kappa)}}
\,|\varphi'(z)|^2\diff A(z) <+\infty,
\end{equation}
where $M$ is given by
$$M=(1+\kappa)(2+\kappa)^{-(2+\kappa)/(1+\kappa)}\const^{1/(\kappa+1)}
|\lambda|^{-1/(\kappa+1)}.$$
A moment's reflection based on the restriction placed on the parameters
$\const$ and $\lambda$ gives that 
\begin{equation}
\label{est2-1}
\int_{\D}
\bigg| \Big[ \frac{z\varphi'(z)}{\varphi(z)}\Big]^{-\tau} \bigg|\,
\big(1-|z|^2\big)^{-\re \tau +R(\tau)}
\,|\varphi'(z)|^2\diff A(z) <+\infty
\end{equation}
holds so long as $R(\tau)$ satisfies
$$R(\tau)>R_0(\tau)=\inf_{0<\kappa<1}
\bigg(\frac{9\e 4^\kappa}{\kappa}\bigg)^{1/(1+\kappa)}
\,\frac{(1+\kappa)\Gamma\big(\frac{1-\kappa}{1+\kappa}\big)}{(2+\kappa)
\Gamma\big(\frac{1}{1+\kappa}\big)^2}\,|\tau|^{(2+\kappa)/(1+\kappa)},$$
where the equality defines $R_0(\tau)$.
We would like to estimate effectively $R_0(\tau)$ as $\tau\to0$.
We realize that if we pick, for small $|\tau|$, 
$$\kappa=\frac{1}{\log\frac1{|\tau|}},$$
then we readily obtain
$$R_0(\tau)\le\bigg[\frac{9\e^2}{2}+o(1)\bigg]\,
|\tau|^2\log\frac{1}{|\tau|}\quad\text{as}\quad |\tau|\to0,$$ 
as claimed.



\begin{thebibliography}{1}

\bibitem{AH} D. R. Adams, L. I. Hedberg, \em Function spaces and potential 
theory. \rm 
Grundlehren der Mathematischen Wissenschaften [Fundamental Principles of 
Mathematical Sciences], {\bf314}. Springer-Verlag, Berlin, 1996. 


\bibitem{carl-jon} L. Carleson, P. W. Jones, \em On coefficient problems for
univalent functions and conformal dimension, \rm Duke Math. J. {\bf 66}
(1992), 169--206.


\bibitem{DV} O. Dragi\v cevi\'c, A. Volberg, {\em Bellman function, 
Littlewood-Paley estimates and asymptotics for the Ahlfors-Beurling operator 
in $L\sp p(\mathbb C)$}. Indiana Univ. Math. J. {\bf54} (2005), 971--995. 

\bibitem{Dur} P. L. Duren, \em Univalent functions. \rm 
Grundlehren der Mathematischen Wissenschaften [Fundamental Principles of 
Mathematical Sciences], {\bf259}. Springer-Verlag, New York, 1983. 


\bibitem{graf} L. Grafakos, \em Classical and modern Fourier analysis. \rm
Pearson Ed., Upper Saddle River, New Jersey, 2004.

\bibitem{Gron} T. H. Gr\"onwall, {\em Some remarks on conformal 
representation.}  Ann. of Math. (2) {\bf 16} (1914-1915), 72--76.

\bibitem{HedShi} H. Hedenmalm, S. Shimorin, {\em Weighted Bergman spaces and
the integral means spectrum of conformal mappings}. Duke Math. J. {\bf127}
(2005), 341--393.

\bibitem{HedShiadd} H. Hedenmalm, S. Shimorin, {\em On the universal 
integral means spectrum of conformal mappings near the origin}.  
Proc. Amer. Math. Soc. {\bf 135}  (2007),  no. 7, 2249--2255.

\bibitem{HedSol} H. Hedenmalm, A. Sola {\em Spectral notions for conformal 
maps: a survey}. Computat. Methods Funct. Theory {\bf8} (2008), No. 2, 
447--474.

\bibitem{JonMak} P. W. Jones, N. G. Makarov, {\em Density properties of 
harmonic measure}, Ann. of Math. (2) {\bf 142}  (1995),  no. 3, 427--455. 

\bibitem{mak} N. G. Makarov, {\em Fine structure of harmonic measure},
St. Petersburg Math. J. {\bf 10}  (1999),  no. 2, 217--268.

\bibitem{Pom} Ch. Pommerenke, {\em Boundary Behaviour of Conformal Maps},
Grundlehren der mathematischen Wissenschaften {\bf299},
Springer-Verlag, Berlin, 1992.

\bibitem{Zyg} A. Zygmund, {\em On certain lemmas of Marcinkiewicz and 
Carleson}.  J. Approximation Theory {\bf 2} (1969), 249--257. 

\end{thebibliography}
\end{document}